\documentclass[11pt,draft]{amsart}
\usepackage{amsmath,amsfonts,amssymb,amscd,amsthm}
\newtheorem{prop}{Proposition}[section]
\newtheorem{thm}[prop]{Theorem}

\newtheorem{lem}[prop]{Lemma}

\theoremstyle{definition}
\newtheorem{de}[prop]{Definition}

\theoremstyle{remark}
\newtheorem*{remark}{Remark}            
\def\dR{\text{\rm dR}}

\def\complex{{\mathbb C}}

\def\PP{{\mathbb{P}}}

\def\real{{\mathbb R}}

\def\CP{{\complex\PP}}

\def\J{{\mathcal{J}}}

\def\cS{{\mathcal{S}}}

\def\cl{\mathop{\rm cl}\nolimits}

\def\dR{\mathop{\rm dR}\nolimits}

\def\Hom{\mathop{\rm Hom}\nolimits}
\def\id{\mathop{\rm id}\nolimits}

\def\ker{\mathop{\rm ker}\nolimits}

\def\rank{\mathop{\rm rank}\nolimits}

\def\span{\mathop{\rm span}\nolimits}

\def\smallnegint{\mathop{\int\mkern-13mu
        \raise.5ex\hbox{${\scriptscriptstyle\diagup}$}}\nolimits}
\begin{document}

\title{Locally holomorphic maps yield symplectic structures}
\author{Robert E. Gompf} \thanks{Partially supported by NSF grants 
DMS-9802533 and DMS-0102922.}
\address{Department of Mathematics, The University of Texas at Austin, 
1 University Station C1200, Austin, TX 78712-0257}
\email{gompf@math.utexas.edu}

\begin{abstract}
For a smooth map $f:X^4\to\Sigma^2$ that is locally modeled by holomorphic 
maps, the domain is shown to admit a symplectic structure that is symplectic 
on some regular fiber, if and only if $f^*[\Sigma]\ne0$. 
If so, the space of symplectic forms on $X$ that are symplectic on all 
fibers is nonempty and contractible.
The cohomology classes of these forms vary with the maximum possible freedom 
on the reducible fibers, subject to the obvious constraints. 
The above results are derived via an analogous theorem for locally 
holomorphic maps $f:X^{2n}\to Y^{2n-2}$ with $Y$ symplectic.
\end{abstract}

\maketitle

\section{Introduction}

Research in the past dozen years has uncovered an intimate relationship 
between the differential topology of closed 4-manifolds and their symplectic 
structures. 
The latter are closed, nondegenerate 2-forms, and have been shown to exist 
on many 4-manifolds (e.g., \cite{G1}). 
Many other 4-manifolds do not 
admit such structures \cite{Ta}, \cite{Sz}, \cite{FS}, even though they 
may be homeomorphic (but not diffeomorphic) to symplectic manifolds.
There has been much recent work aimed at understanding which 
4-manifolds admit symplectic structures, as well as the range of values of 
the Chern class $c_1(\omega)$ and cohomology class $[\omega]$ of symplectic 
forms $\omega$ on a fixed manifold (e.g., \cite{MT}, \cite{S}, \cite{LL}). 
For example, a closed 4-manifold admits a symplectic structure if and only
if it admits a fibrationlike structure called a {\em Lefschetz pencil\/} 
\cite{D}, \cite{G2}. 
In fact, a Lefschetz pencil (with fibers suitably intersecting the base locus)
determines a symplectic structure up to isotopy, 
and a dense subset of all symplectic forms is realized this way up to scale.
In this article, we use tools from \cite{G2} to investigate a related problem:
We show that any smooth $f:X^4\to \Sigma^2$ that is locally modeled by 
holomorphic maps allows us to construct symplectic forms on $X$, provided 
that $f^* [\Sigma]\ne0$. 
We show that the space of symplectic forms suitably compatible with $f$ is  
contractible, but that there is much freedom in the class $[\omega]$ when 
$f$ is sufficiently singular. 
We also investigate the corresponding problem in higher dimensions 
(Theorem~\ref{thm:5.2}). 

The maps of interest are defined as follows. 

\begin{de} \label{def:5.1}
A map $f:X^{2n} \to Y^{2m}$ between smooth, oriented manifolds is
{\em locally holomorphic\/} if for each $x \in X$ there are smooth,
orientation-preserving coordinate charts about $x$ and $f(x)$ (into
$\complex^n$ and $\complex^m$, respectively) in which $f$ is given by a
holomorphic map. 
\end{de}

\noindent
The fibers $f^{-1}(y)$ of a locally holomorphic map $f$ 
are homeomorphic to CW-complexes.
If each component of a compact fiber has (real) dimension~2, 
then the fiber is the image of some closed 
surface $F$ under a map that is a smooth embedding outside a finite
subset of $F$ whose image we will denote by $K_y \subset X$. 
Near $K_y$, $f^{-1}(y)$ is conelike.  
As in algebraic geometry, we will call the
image of each connected component of $F$ an {\em irreducible
component\/} of $f^{-1}(y)$.  
Note that some irreducible components
may have multiplicities $>$ 1, so that they are contained in the
critical set of $f$.
In our case of primary interest ($X$ compact, $n=2$, $m=1$), a 
(nonconstant) locally holomorphic map is locally modeled by $f(z,w)=z^d$ 
where $d\ge1$ is the (local) multiplicity of the fiber $f^{-1}(0)$, 
except on the finite set $K=\bigcup K_y$. 
Thus, the tangent spaces to the fibers form a real 2-plane bundle $L$ 
over $X-K$, and this is canonically oriented (by the preimage orientation 
at regular points). 
There are simply connected 4-manifolds that admit no complex structure, 
but do admit locally holomorphic maps to $S^2$ (e.g., \cite{GS}), even 
if we require all critical points to be of the simplest type (quadratic). 
One can construct more complicated examples, starting from any finite  
collection of (connected) singular fibers of holomorphic maps of a 
fixed generic fiber genus, by gluing together their tubular neighborhoods,  
and extending the resulting singular fibration from $D^2$ to $S^2$ by 
adding quadratic critical points to cancel the geometric monodromy. 
Note that each original singular fiber could have many irreducible 
components, with differing multiplicities. 

We can now begin to state our main results. 

\begin{thm}\label{thm:5.1} 
Let $X$ be a closed, oriented, connected 4-manifold with a locally 
holomorphic map $f:X\to \Sigma$ to a closed, oriented surface. 
Then $X$ admits a symplectic structure that is symplectic on some 
regular fiber, if and only if $f^* [\Sigma] \ne0\in H_{\dR}^2 (X)$. 
\end{thm}

\noindent 
The hypothesis $f^*[\Sigma] \ne 0$ is equivalent to requiring a
generic fiber to be nontrivial in $H_2(X;\real)$. 
This is clearly necessary for the existence of symplectic structures on
$X$ as above, since 
$\omega$ being symplectic on a generic fiber $f^{-1}(y)$ implies
$\langle [\omega], [f^{-1}(y)] \rangle \ne 0$.
(Sufficiency follows from Theorem~\ref{contractible}  below.)  
However, the hypothesis $f^*[\Sigma]
\ne 0$ is automatically satisfied whenever the generic fiber is not a
torus (or disjoint union of tori).  
To see this, note that the above bundle $L$ over $X - K$ defines 
an Euler class $e(L) \in H^2(X-K) \cong H^2(X)$. 
For a generic fiber, $\langle e(L), f^{-1}(y) \rangle = \langle
e(Tf^{-1}(y)), f^{-1}(y) \rangle = \chi(f^{-1}(y)) \ne 0$ unless
$f^{-1}(y)$ is a union of tori 
(since its components must be diffeomorphic).  
If we also require each
fiber to have a neighborhood on which $f$ is globally modeled by a
holomorphic map, then the classification of elliptic fibrations
implies that $f^* [\Sigma] \ne 0$ unless $f$ is (up to blowups) a
Seifert fibration (i.e., made from an honest torus bundle by adding
smooth multiple fibers, and possibly composing with a branched
covering map of $\Sigma$ to allow disconnected fibers). 

We would like to study the space of symplectic structures on $X$ that are 
symplectic on all fibers of $f$. 
However, this condition makes no sense on the finite set $K$ where the 
fibers are singular. 
Instead, we use the following  proposition, which is proved at 
the end of this section. 

\begin{prop}\label{Jstar} 
A locally holomorphic map $f:X^4\to\Sigma^2$ canonically determines a 
complex bundle structure $J^*$ on $TX| K$. 
This is obtained by restricting the complex structures inherited from 
any choices of charts as in Definition~\ref{def:5.1}, or more generally, 
by restricting any $C^0$ almost-complex structure $J$, defined on a 
neighborhood $U$ of $K$, for which the fibers of $f|U-K$ are $J$-holomorphic.
\end{prop}

\noindent
Recall that an almost-complex structure $J$ on $U$ 
is defined to be a complex structure on the tangent bundle 
$TU$, i.e., a bundle map covering $\id_U$ with $J\circ J= -\id_{TU}$. 
We always assume almost-complex and symplectic structures respect the given 
orientations. 
Thus, in the proposition, 
$J$ determines the given orientation on $U$, and the tangent spaces 
to the regular fibers are $J$-complex lines in the preimage orientation. 
We only need almost-complex structures 
for a crude level of directional control in the tangent spaces. 
Thus, it is convenient to ignore differentiability and only require 
continuity of $J$, although we indicate where 
smoothness can be imposed for stronger conclusions.

\begin{thm}\label{contractible}
For $f:X\to\Sigma$ as in Theorem~\ref{thm:5.1} 
with $f^* [\Sigma] \ne0$, let $\cS$ be the space of symplectic forms on $X$ 
that are symplectic on the (canonically oriented) fibers of $f|X-K$ and 
tame $J^*$ on $K$. 
Then: 
\begin{itemize}
\item[a)] $\cS$  is nonempty and contractible.
\item[b)] $\cS$ is also characterized as being the space of symplectic forms 
taming global $C^0$  (or $C^\infty$)
almost-complex structures $J$ on $X$ such that each 
regular fiber of $f$ (and hence each fiber of $f|X-K$) is $J$-holomorphic. 
\end{itemize}
\end{thm}

\noindent
For taming, see Definition~\ref{def:2.1}. 
There  is considerable freedom to choose the topology on $\cS$. 
For example, we can choose any $C^k$-topology or Sobolev topology in between. 
See \cite[Theorem 2.11(b) and subsequent discussion]{G2} for further details. 
Note that Theorem~\ref{contractible} immediately implies 
Theorem~\ref{thm:5.1}.  

To understand our freedom to choose the cohomology class $[\omega]\in 
H_{\dR}^2(X)$, first note some obvious constraints: 
If $F_1,\ldots,F_n$ are the (canonically oriented) irreducible components 
of a single connected component of a fiber, and $F_i$  has multiplicity 
$m_i$, then (as we will see) $\sum m_i [F_i]\in H_2(X;\real)$ must be a 
rational multiple of the generic fiber class $[f^{-1}(y)]$ (where the 
multiplier $q$ is 1 if the fibers are connected). 
Thus we must have $\sum m_i\langle\omega,F_i\rangle = q\langle \omega,
f^{-1}(y)\rangle$. 
Furthermore, for $\omega\in \cS$ each $F_i-K$ must be symplectic, so 
$\langle \omega,F_i\rangle >0$. 
These turn out to be the only constraints on the areas 
$\langle \omega,F_i\rangle$. 

\begin{thm}\label{areas}
For fixed $J$ as in Theorem~\ref{contractible}(b), there is a form 
$\omega\in\cS$ taming $J$, such that the areas 
$\langle\omega,F_i\rangle$ of irreducible components of fibers of $f$ 
take any preassigned values subject to the above constraints.
\end{thm}

\noindent
That is, the (suitably weighted) areas of the irreducible components of 
each  connected component of each fiber can be distributed in any 
preassigned manner. 

Since $\cS$ is connected, all $\omega\in\cS$ have the same Chern class 
$c_1(\omega)$. 
We can compute this class using some $J$ as in Theorem~\ref{contractible}(b) 
and the $J$-complex line bundle $L$ over $X-K$ of tangent planes to fibers. 
If the critical set of $f$ in $X$ is finite, then 
by definition $c_1(\omega) = c_1(TX, J) = c_1(L \oplus f^*T\Sigma) = e(L) +
\chi(\Sigma)f^*[\Sigma] = e(L) + \chi(\Sigma)PD[f^{-1}(y)]$
for a generic fiber $f^{-1}(y)$ (where PD is Poincar\'e duality). 
The general case is obtained by adding a term $(1-m_i)PD[F_i]$ for each 
irreducible component $F_i$ of multiplicity $m_i$ in each singular fiber.
(A vector field on $\Sigma$ which is nonzero near the critical values 
of $f$ lifts to one on $L^\perp$  with index $1-m_i$ near $F_i$, since 
sources (index~1) at critical values lift to sources.) 

We prove Theorems~\ref{contractible} and \ref{areas} in Section~3, 
using tools from \cite{G2}. 
Since these tools are applicable in 
arbitrary dimensions, we proceed by first proving a theorem that  holds
in all dimensions. 
This theorem says that under suitable hypotheses, 
a locally holomorphic map to a 
symplectic manifold, with 2-dimensional fibers, determines a deformation 
class of symplectic structures on its domain (a {\em deformation\/} being
a smooth family $\omega_t$, $0\le t\le 1$, of symplectic forms). 
The precise statement (Theorem~\ref{thm:5.2}) requires further definitions, 
and is the subject of the next section. 

  
\begin{proof}[Proof of Proposition~\ref{Jstar}] 
For a locally holomorphic map $f:X^4\to\Sigma^2$ and $x\in K$ a singular 
point of a fiber, we must show that all almost-complex structures defined 
near $K$, making the fibers of $f$ $J$-holomorphic, have the same 
restriction $J^*$ to $T_xX$. 
Identify a 4-ball neighborhood $U$ of $x$ in $X$ with a neighborhood 
of 0 in $\complex^2$ as in Definition~\ref{def:5.1}. 
It is not hard to show that a linear complex structure on $\real^{2n}$ 
($n\ne1$) is determined by its 1-dimensional (oriented) complex subspaces 
\cite[Lemma 4.4(a)]{G2}, so it suffices to show that every 
complex line at $0$ in $\complex^2$ can be written as $\lim \ker df_{x_i}$ 
for some sequence of regular points $x_i\to0$. 
We can assume $U \cap K = \{x\}$, so the bundle $L$ of tangents to fibers
of $f$ is defined on $U-\{0\} \subset \complex^2$ and determines a
holomorphic map $\varphi:U - \{0\} \to \CP^1$, whose homogeneous
coordinates $\varphi_i$, $i=1,2$, are obtained from 
$\frac{\partial f}{\partial z_i}$ by removing common factors (to remove the
singularities along the critical set).
Up to homotopy, continuous maps $U-\{0\} \simeq S^3 \to \CP^1$ are
classified by $\pi_3(S^2) \cong \mathbb{Z}$, and we can compute this integer
invariant $h(\varphi)$ by the Thom-Pontrjagin construction: The fibers
$\varphi_i^{-1}(0)$ of $\varphi$ each intersect a small $S^3 \subset U
- \{0\}$ in an oriented link $L_i$ (possibly with multiplicities), and
$h(\varphi) = \ell k(L_1, L_2)$ is their linking number.
If $0 \not \in \varphi_i^{-1}(0)$ for some $i$ then $\varphi$
determines a holomorphic map $U-\{0\} \to \complex$ (for $U$ sufficiently 
small), so $L$ extends holomorphically over $0$ by Hartogs' Theorem.
But $L|U-\{0\}$ is tangent to the fibers of $f$, so $L$ is integrable
on $U$ (since the Frobenius condition is closed), and $f^{-1}(0)$ is a
smooth leaf of the foliation, contradicting our assumption that $x \in K$.
Thus, $0 \in \varphi_1^{-1}(0) \cap \varphi_2^{-1}(0)$, so $h(\varphi)
= \ell k(L_1, L_2) = \varphi_1^{-1}(0) \cdot \varphi_2^{-1}(0) > 0$.
In particular, $\varphi$ is surjective on arbitrarily small 3-spheres about 0, 
and so the required sequence $(x_i)$ is easy to construct.
\end{proof}

\section{Arbitrary dimensions}

We begin with some terminology for relating
symplectic and complex structures.

\begin{de} \label{def:2.1}\cite{G2}
Let $T:V \to W$ be a linear transformation between finite-dimensional
real vector spaces, and let $\omega$ be a skew-symmetric bilinear form on $W$.  
A linear  complex structure $J:V \to V$ will be called 
$(\omega,T)$-{\em tame\/} if $T^*\omega(v, Jv) > 0$ for all $v \in V
- \ker T$.  If, in addition, $T^*\omega$ is $J$-invariant (i.e., 
$T^*\omega(Jv, Jw) = T^*\omega(v,w)$ for all $v,w \in V$), we
will call $J$ $(\omega, T)$-{\em compatible\/}.  For a $C^1$-map $f:X
\to Y$ between manifolds, with a 2-form $\omega$ on $Y$, an
almost-complex structure $J$ on $X$ will be called
$(\omega,f)$-{\em tame\/} (resp. $(\omega, f)$-{\em compatible\/}) if it
is $(\omega,df_x)$-{\em tame\/} (resp. $(\omega,df_x)$-{\em compatible\/}) 
for each $x \in X$.  If $T = \id_V$ or $f = \id_X$, we will shorten
the terminology to $\omega$-{\em tame\/} and $\omega$-{\em compatible\/}. 
\end{de}

\noindent 
The last sentence of the definition is standard terminology.  
If $\omega$ tames some $J$  (i.e., $T=\id_V$ and 
$J$ is $\omega$-tame), then $\omega$ is positive on all $J$-complex lines,
and it is  obviously nondegenerate (that is, every $v\ne0$ 
pairs nontrivially with something), 
so a closed, taming 2-form $\omega$ is automatically symplectic.  
Both the taming and compatibility
conditions are preserved under taking convex combinations $\sum t_i
\omega_i$  (all $t_i \geq 0, \sum t_i = 1$) for fixed $J$. 
However, taming is more natural for our purposes than compatibility, since 
$\omega$-taming is preserved under small perturbations of $\omega$ and $J$. 
It is an open question whether compatibility can be replaced by taming 
throughout the paper. 
(It can be done, for example, if Question~4.3 of \cite{G2} has an 
affirmative answer.) 
In the special case where $\dim_{\real}W=2$ and $\omega$ respects a preassigned 
orientation on $W$, then $\omega$ is determined up to positive scale, so 
$(\omega,T)$-taming and $(\omega,T)$-compatibility are equivalent, 
independent of choice of $\omega$, and equivalent to the condition that 
$\ker T$ (in the preimage orientation if $T\ne0$, where $J$ orients $V$)
be a $J$-complex subspace of $V$. 
In particular: 

\begin{prop}\label{tame}
For $f:X\to\Sigma$ as in Theorem~\ref{contractible}, $J$ satisfies the 
condition of (b) (every regular fiber is $J$-holomorphic) if and only if $J$ 
is $(\omega_\Sigma,f)$-tame (or $(\omega_\Sigma ,f)$-compatible), where 
$\omega_\Sigma$ is any (positive) area form on $\Sigma$.\qed
\end{prop}

\noindent
Such structures $J$ are easy to construct: 
Split $T(X-K)$ as $L\oplus L'$ where $L$ is tangent to the fibers and $L'$
is complex near $K$ for some locally defined complex structure as in 
Definition~\ref{def:5.1}, then declare 
$L$ and $L'$ to be complex line bundles. 
However, it is more difficult to arrange $J$ to tame a preassigned 
$\omega \in \cS$. 
Thus, while it is easy to see that any $\omega$ as in 
Theorem~\ref{contractible}(b) lies in $\cS$, the converse takes more work.

In higher dimensions, we need to strengthen our local holomorphicity 
condition:

\begin{de}\label{holomorphic} 
For a symplectic form $\omega$ on $Y$, a locally holomorphic map $f:X\to Y$ 
is called {\em $\omega$-compatibly locally holomorphic\/} if the charts in 
Definition~\ref{def:5.1}  
can be chosen so that $\omega$ is compatible with the standard complex 
structure on $\complex^m$.
\end{de}

\noindent
This is automatically satisfied when $\dim_{\real} Y=2$, or for 
holomorphic charts on a K\"ahler $(Y,\omega)$. 
It implies that the corresponding local 
complex structures on $X$ are $(\omega,f)$-compatible. 

We also need a technical condition from \cite{G2} to control the 
behavior of fibers near critical points; this is again vacuous for 
nonconstant locally holomorphic maps $f:X^4\to \Sigma^2$. 
Suppose $E,F \to X$ are real vector bundles over a metrizable 
topological space,
and $T:E \to F$ is a  section of the bundle $\Hom(E,F)$. 
In our main application, these will be induced by a $C^1$-map $f:X \to Y$
between manifolds, with $T= df:TX \to f^*TY$. Motivated by
this example, we call a point $x \in X$ {\em regular\/} if $T_x:E_x \to
F_x$ is onto and {\em critical\/} otherwise.  
Let $P \subset E$ be the
closure $\cl( \bigcup \ker T_x)$, where $x$ varies over all the regular
points of $T$ in $X$, and let $P_x = P \cap E_x$.  Thus, $P_x = \ker
T_x$ if $x$ is regular, and otherwise $P_x \subset \ker T_x$ consists
of limits of sequences of vectors annihilated by $T$ at regular
points. 

\begin{de} \label{def:2.2}\cite[Definition 2.2]{G2}
A point $x \in X$ is {\em wrapped\/} if $\span P_x$ has (real)
codimension at most 2 in $\ker T_x$.
\end{de}

\noindent
In our application, this condition automatically holds away from fiber
components with multiplicities: 

\begin{prop} \label{prop:2.3}\cite[Proposition 2.3]{G2}
Suppose that in a neighborhood of a critical point $x \in X$, $T$ is
given by $df$, for some holomorphic map $f:U \to
\complex^{n-1}$ with $U$ open in $\complex^n$.  
If each fiber
$f^{-1}(y)$ intersects the critical set  of $f$ in at most a finite
set, then $x$ is wrapped.  
In fact, $P_x = \ker T_x$.\qed  
\end{prop}

\noindent
Note  that the
proposition becomes false without the finiteness hypothesis, e.g.,
$n=3$, $f(x,y,z)=(x^2, y^2)$ at $(0,0,0)$.  
For $n=2$, $P_x = \ker T_x$ unless $f$ is constant or $x$ is a smooth
point of a fiber component 
with multiplicity $>1$, cf. proof of Proposition~\ref{Jstar},
but every point of a nonconstant locally holomorphic map $f:X^4\to\Sigma^2$ 
will be wrapped (since regular points are dense, implying 
$\dim_{\real}\text{span}\, P_x \ge2$ everywhere). 

To state the main theorem, we must first orient the fibers. 

\begin{lem}\label{fibers}
Let $f:X^{2n}\to Y^{2n-2}$ be a locally holomorphic map, all of whose fibers  
have real dimension~$2$. 
Fix $y\in Y$. 
\begin{itemize}
\item[a)] 
For each $x\in f^{-1}(y)-K_y$, there is a sequence $x_i\to x$ of regular
points in $X$ for which $T_xf^{-1}(y) = \lim\ker df_{x_i}$. 
\item[b)] 
The surface $f^{-1}(y)-K_y$ is canonically oriented, by sequences as in (a) 
and the preimage orientation on each $\ker df_{x_i}$.
\end{itemize}
\end{lem}

\begin{proof} 
For (a) it suffices to find such a sequence $(x_i)$ for each $x$ in a dense 
subset of $f^{-1}(y)$. 
After restricting to suitable neighborhoods, we may assume $f$ is holomorphic. 
The critical values of $f$ lie in a local subvariety of $Y$ with positive 
codimension, so there is a holomorphic disk $D$ in $Y$ 
(with $\dim_{\complex} D=1$) centered at $y$ and disjoint from 
the critical set elsewhere. 
After suitably shrinking $D$, we conclude that $D\subset f(X)$. 
(Otherwise $f(X)$ lies in a variety of positive codimension, so $f$ has 
fibers of real dimension $>2$.) 
Thus, $f^{-1} (D)$ is locally a variety of complex dimension~2. 
After we resolve the singularities of $f^{-1}(D)$, the required sequences 
obviously exist on the resulting complex surface, since generically $f$ is 
given locally by $f(z,w)=z^d$ with fibers parallel to the $w$-axis. 
The sequences push down to $f^{-1}(D)$ as required. 
(This follows a suggestion of S.~Keel.)

To prove (b), note that the preimage orientation on each $\ker df_{x_i}$ 
is also its complex orientation (for the complex structure determined by 
our choice of charts in Definition~\ref{def:5.1}). 
These then limit to the complex orientation on $T_xf^{-1}(y)$, which is 
independent of the choice of $(x_i)$ converging to $x$. 
Since this orientation is now determined by a fixed $x_i\to x$, it is also 
independent of our choice of charts, and hence global and canonical.
\end{proof}  

\begin{thm} \label{thm:5.2}
Let $X^{2n}$, $Y^{2n-2}$ be closed, oriented manifolds with a symplectic
form $\omega_Y$ on $Y$ and an $\omega_Y$-compatibly locally
holomorphic map $f:X \to Y$, all of whose fibers have real dimension~2.  
Suppose there is a class $c \in
H^2_{\dR}(X)$ evaluating positively on each irreducible component of each 
fiber of $f$ (canonically oriented). 
If $n \geq 3$, assume that all critical points of $f$ are wrapped. 
Then: 
\begin{itemize}
\item[a)] $X$ admits a symplectic structure.  
In fact, there is a unique deformation class of
symplectic forms on $X$  containing representatives that tame 
$(\omega_Y, f)$-compatible, global $C^0$ 
almost-complex structures $J$.   
This still holds if we require $J$ to be $C^\infty$, with a $C^\infty$, 
$\omega_Y$-compatible structure on $f^* TY$ making $df$ $J$-complex.
\item[b)] For any fixed $C^0$, $(\omega_Y,f)$-tame $J$ on $X$, 
the convex open cone in $H_{\dR}^2(X)$ consisting of classes of 
symplectic forms taming $J$ contains all classes $tc+f^*[\omega_Y]$ 
for $c$ as above and $t>0$ sufficiently small (depending on $c$).  
\end{itemize}
\end{thm}

%
%
%

\begin{remark}
When $n=3$, the condition of wrapped critical points can be dropped, and 
existence and (b) still follow. 
In fact the entire theorem still holds then, if we also require the 
structures $J$ in (a) to be $(\omega_Y,df)$-extendible as in \cite{G2}, along 
preassigned sequences converging to the unwrapped critical points. 
The proof is identical to that of Theorem~2.7, once we augment 
Lemma~3.2 by \cite[Addendum 3.3]{G2} and note that the charts of 
Definition~2.3 are automatically $(\omega_Y,df)$-extendible.
The resulting deformation class in (a) is independent of the choice of 
sequences, since any $J$ satisfying all the conditions of (a) is 
$(\omega_Y,df)$-extendible for all sequences. 
\end{remark}

\section{Proofs}

It remains to prove 
Theorems~\ref{thm:5.2}, \ref{contractible} and \ref{areas}. 
We need two results from \cite{G2}. 
The first allows us to construct symplectic structures on the domain $X$ 
of a map to a symplectic manifold, in the presence of a suitable 
almost-complex structure $J$ on $X$. 
The second allows us to construct $J$. 

\begin{thm} \label{thm:3.1}
\cite{G2}
Let $f:X \to Y$ be a smooth map between closed manifolds. 
Suppose that $\omega_Y$ is a symplectic form on $Y$, and $J$ is a
continuous, $(\omega_Y, f)$-tame almost-complex structure on $X$.  
Fix a class $c \in H^2_{\dR}(X)$.  
Suppose that for each $y \in Y$,
$f^{-1}(y)$ has a neighborhood $W_y$ in $X$ 
with a closed 2-form $\eta_y$  such that $[\eta_y] = c|W_y
\in H^2_{\dR}(W_y)$ and such that $\eta_y$ tames $J$ on each
of the complex subspaces $\ker df_x$, $x \in W_y$. 
Then there is a closed 2-form $\eta$ on $X$ with $[\eta] = c \in
H^2_{\dR}(X)$, and such that for all sufficiently small $t>0$
the form $\omega_t = t\eta + f^*\omega_Y$ on $X$ tames $J$ (and hence is
symplectic).\qed
\end{thm}

This is \cite[Theorem 3.1]{G2}, restricted to the case with $C=\emptyset$. 
The main idea of the proof goes back to Thurston \cite{T} in the case 
of surface bundles. 
We cannot directly splice the forms $\eta_y$ by a partition of unity, without 
losing closure of the forms. 
Instead, we subtract off a global representative $\zeta$ of $c$ to obtain exact 
forms, then splice via the corresponding 1-forms and add $\zeta$ back in. 
The resulting closed 2-form $\eta$ is nondegenerate on each $\ker df_x$ by 
convexity of the $J$-taming condition. 
(This use of $J$ to control nondegeneracy is the innovation allowing us to 
deal with general critical points of $f$.) 
The form $f^*\omega_Y$ provides nondegeneracy for $\omega_t$ in the 
remaining directions. 

Splicing together local almost-complex structures is harder. 
Here we only state a simplified version of 
\cite[Lemma 3.2]{G2} without proof. 
As preceding Definition~\ref{def:2.2}, we let $E, F \to X$ be
real vector bundles over a metrizable space, with fiber dimensions
$2n$ and $2n-2$ respectively, and this time equipped with fiber orientations.  
We again fix a section $T:E \to F$ of $\Hom (E,F)$.  
A {\em 2-form\/} on $E$ or $F$ means  a continuous choice of 
skew-symmetric bilinear forms on the fibers.

\begin{lem} \label{lem:3.2}
\cite{G2}
For $E,F \to X$ and $T$ as above with $X$ compact, 
suppose that the regular points of $T$ are dense in $X$, 
and let $\omega_F$ be a nondegenerate 2-form on $F$ (inducing
the given fiber orientation).  
Suppose that each $x \in X$ has a
neighborhood $W_x$ with an $(\omega_F, T)$-compatible complex (bundle) 
structure on $E|W_x$. 
\begin{itemize}
\item[a)] If $n \geq 3$, assume each critical point of $T$ is wrapped. 
Then the space $\J$ of $(\omega_F, T)$-compatible complex structures  on $E$ 
is nonempty and contractible (in the $C^0$-topology). 
For any $(\omega_F,T)$-compatible structure $J_C$ defined near a closed 
subset  $C\subset X$, there are elements of $\J$ agreeing with $J_C$ on $E|C$.
\item[b)] Fix a 2-form $\omega_E$ on $E$ and a complex structure on $F$. 
Then (a) remains true if each complex structure  on $E$ and its restrictions 
is required to be $\omega_E$-tame and to make $T$ complex linear.\qed 
\end{itemize}
\end{lem}

\begin{proof}[Proof of Theorem~\ref{thm:5.2}] 
For $f:X\to Y$ as  given, Lemma~\ref{lem:3.2}(a) (with $T= df: TX\to 
f^* TY$) implies that the space $\J$ of $(\omega_Y,f)$-compatible 
almost-complex structures on $X$ is nonempty and contractible. 
To see this, note that regular points of $f$ are dense (e.g., by 
Lemma~\ref{fibers}(a)), and that the local complex structures on $X$ given by 
Definition~\ref{holomorphic} and subsequent text are automatically 
$(\omega_Y,f)$-compatible. 
The relevant part of the existence proof of Lemma~3.2(a) is smooth, 
and automatically produces elements of $\J$ satisfying all the 
conditions of Theorem~\ref{thm:5.2}(a).
Now we check that each (canonically oriented) fiber $f^{-1}(y)-K_y$ is 
a $J$-holomorphic curve for each $J\in \J$, or more generally, for any 
$(\omega_Y,f)$-tame (positively oriented) $C^0$ almost-complex structure $J$ 
on $X$. 
This follows from Lemma~\ref{fibers} once we verify that each 
$\ker df_x$, with $x$ regular, is a $J$-complex line in the preimage 
orientation. 
But Definition~\ref{def:2.1} implies $\ker df_x$ is $J$-invariant, so 
$T_{f(x)}Y \cong T_xX/\ker df_x$ inherits an $\omega_Y$-tame complex 
structure $f_*J$. 
After homotoping $\omega_Y$ through taming structures to one that is 
compatible with $f_*J$, it is easy to verify that $f_*J$ induces the 
same orientation on $T_{f(x)} Y$ as $\omega_Y$. 
Thus $J$ induces the preimage orientation on $\ker df_x$.

For any fixed  $(\omega_Y,f)$-tame $J$ on $X$, 
such as any $J\in \J$, 
we wish to apply Theorem~\ref{thm:3.1} with $c$ as given.  
The argument follows the method of \cite[Theorem 2.11(b)]{G2}, 
but with complications arising from fibers with irreducible components 
lying in the critical set of $f$.  
For $y\in Y$, 
let $K_y' \supset K_y$ be a finite subset of $f^{-1}(y)$
intersecting each irreducible component nontrivially, and let $\sigma$
be a closed $2$-form taming $J$ on some neighborhood $W$ of $K_y'$ in $X$. 
Since $f^{-1}(y)-K_y$ is $J$-holomorphic, 
$\sigma|(f^{-1}(y) - K_y) \cap W$ is a positive area form,
and we can extend this to an area form on $f^{-1}(y) - K_y$.  
We can easily arrange $\int_{F_i} \sigma = \langle c, F_i\rangle
> 0$ for each irreducible component $F_i$ of $f^{-1}(y)$.  
Let $F^* \subset f^{-1}(y)$ be a compact surface with boundary, obtained
by deleting a neighborhood of $K_y'$ whose closure lies in $W$. 
Since $F^*$ is $J$-holomorphic, $TF^* \subset TX|F^*$ has a complementary
complex subbundle $\nu F^*$.  
This bundle is trivial since $F^*$ has
no closed components, so we can use it to identify a tubular
neighborhood $N$ of $F^*$ with $F^* \times \complex^{n-1}$, by a map
that is $J$-holomorphic on $TN|F^*$.  
After shrinking $W$, we may assume $N\cap W$ corresponds to 
$(F^*\cap W)\times \complex^{n-1}$. 
The product form $\tau =
\pi_1^* (\sigma|F^*) + \pi_2^* \omega_{\complex^{n-1}}$ on $N$
tames $J$ on $TN|F^*$.  
The form $\sigma - \tau$ on $N \cap W$ is
closed, and it vanishes on $F^*\cap W$, so it is exact. 
Choose a 1-form $\alpha$ on $N \cap W$ with $d\alpha = \sigma - \tau$.  
Subtracting the closed form $\pi_1^*(\alpha |F^* \cap W)$
from $\alpha$ if necessary, we can arrange $\alpha |F^* \cap W = 0$. 
Now the 1-form $\alpha:T(N\cap W) \to \real$ restricts to a 
fiberwise-linear map $\nu F^*|F^* \cap W \to
\real$; by our identification of this bundle with $N \cap W$, we
obtain a smooth map $\varphi : N \cap W \to \real$ with $\varphi|F^*
\cap W = 0$ and $d\varphi = \alpha$ on $TX|F^* \cap W$.  
Subtracting
$d\varphi$ from $\alpha$, we arrange $\alpha =0$ on $TX|F^* \cap W$.
Choose a map $\rho:F^* \to [0,1]$ with $\rho = 1$ near $\partial F^*$
and $\rho = 0$ outside $W$, and let $\eta_y = \tau + d((\rho \circ
\pi_1)\alpha)$ on $N$.  
Then the closed form $\eta_y$ agrees with
$\tau$ outside $W$ and extends as $\sigma$ near $f^{-1}(y) - F^*$.
Furthermore, $\eta_y = \sigma$ on $f^{-1}(y)-K_y$. 
On $TX|F^* \cap W$,
$d(\rho \circ \pi_1) \wedge \alpha = 0$ so $\eta_y = \tau + \rho
d\alpha$ is a convex combination of the $J$-taming forms $\tau$ and $\sigma$. 
Thus,  $\eta_y$ tames $J$ on $TX|f^{-1}(y)$ and hence on a regular 
neighborhood $W_y$ of $f^{-1}(y)$ in $X$.  
Since $\int_{F_i}\eta_y = \int_{F_i}\sigma = \langle c,F_i\rangle$ for 
each $F_i$, so $[\eta_y] = c|W_y$, 
Theorem~\ref{thm:3.1}
gives us a global closed 2-form $\eta$ on $X$ with $[\eta] = c$ and
$\omega = t \eta + f^* \omega_Y$ taming $J$ (hence, symplectic) for
any sufficiently small $t > 0$. 
As required, $[\omega] = tc+f^* [\omega_Y]$, 
so we have proved (b) and the existence part of (a).

To prove uniqueness in (a), we must find a deformation between any
preassigned pair $\omega_s$, $s=0,1$, of forms taming structures  $J_s\in \J$. 
Since $\J$ is contractible, 
we can extend to a family $J_s\in\J$, $0 \leq s \leq 1$.  
For each $s \in (0,1)$, Part~(b) 
(with any convenient choice of $c$) yields a symplectic
form $\omega_s$ taming $J_s$ (not necessarily continuous in $s$).  
Since the taming condition is open and
$X$ is compact, each $\omega_s$, $s \in [0,1]$, tames $J_t$ for $t$ in
some neighborhood of $s$.  
Convexity of the taming and closure conditions now allows us to 
smooth the family $\omega_s$
by a partition of unity on $[0,1]$, to obtain the required
deformation of symplectic forms. 
\end{proof}
  
\begin{proof}[Proof of Theorem~\ref{areas}] 
For $f:X\to\Sigma$ and $J$ as in Theorem~\ref{contractible}, we wish to 
construct symplectic structures using Theorem~\ref{thm:5.2}(b). 
We have already observed (following Definition~\ref{holomorphic}) that 
$f$ is $\omega_\Sigma$-compatibly locally holomorphic for any 
$\omega_\Sigma$, and (Proposition~\ref{tame}) $J$ is 
$(\omega_\Sigma,f)$-tame.  
Since $f^*[\Sigma] \neq 0$, $f$ is surjective (onto a component of 
$\Sigma$), so local holomorphicity implies each fiber is 2-dimensional.  
To construct a suitable class $c$, we invoke the
following lemma, whose proof is given at the end of the paper.

\begin{lem} \label{lem:5.4}
For an $n \times n$ real, symmetric matrix $A = [a_{ij}]$, let $G_A$
denote the graph with $n$ vertices $v_1, \ldots, v_n$, and an edge
between any two distinct vertices $v_i, v_j$ whenever $a_{ij} \ne 0$.  
Suppose that (a) $G_A$ is connected, (b) $a_{ij} \geq 0$ whenever
$i \neq j$, and (c) there are positive real numbers $m_1, \ldots, m_n$
such that $\sum_{i=1}^{n} m_i a_{ij} \leq 0$ for all $j$.  
Fix a choice of such numbers $m_i$. 
Then the hypothesis (d), that  the inequality in (c) is strict for some 
$j$, implies $\rank A = n$.
If (d) is not satisfied, then $\rank A = n-1$.
\end{lem}

\noindent
Let $F_1, \ldots, F_n$ denote the irreducible components comprising a
single connected component of a singular fiber $f^{-1}(y_0)$. 
For a disk
$D$ about $y_0$ containing no other critical values, and $y \ne y_0$
in $D$, let $\varphi \in H_2(X;\real)$ be the homology class of the
union of all components of $f^{-1}(y)$ lying in the same component of
$f^{-1}(D)$ as $\bigcup_{i=1}^n F_i$. 
Then $\varphi = \sum_{i=1}^n m_i
[F_i] $, where $m_i > 0$ is the multiplicity of $F_i$.  
Furthermore, $\varphi$ is a positive rational multiple of $[f^{-1}(y)] = PD
f^*[\Sigma]$. 
(Since $f$ restricts to a fiber bundle with connected
total space away from the finite set of critical values in $\Sigma$,
any two components of generic fibers are isotopic.) 
If $A = [a_{ij}]$
is the symmetric $n \times n$ matrix for which $a_{ij}$ is the
intersection number $F_i \cdot F_j$, then $\sum_{i=1}^n m_i a_{ij} =
\varphi \cdot F_j = 0$ for all $j$ (since $f^{-1}(y) \cap F_j = \emptyset$).  
Thus $A$ satisfies (a)--(c) but not (d) of the above
lemma, so $\rank A = n - 1$.  
For any $s_1, \ldots, s_n \in \real$ with
$\sum m_i s_i = 0$, we can now find constants $r_1, \ldots, r_n \in
\real$ such that $\psi  = \sum r_i [F_i]$ satisfies $\psi \cdot F_j =
s_j$ for all $j$.  
Clearly, $\psi$ pairs trivially with every other
irreducible component of every fiber.  
Since $[f^{-1}(y)] \neq 0$ by
hypothesis, we can find a class $c_0$ in $H^2_{\dR}(X)$ with
$\langle c_0, f^{-1}(y) \rangle > 0$ (so $c_0$ is positive on each
component of each generic fiber).  
After adding the Poincar\'e dual of
a suitable class $\psi$ for each connected component of each singular
fiber, we obtain a class $c$ realizing any preassigned values on
irreducible components of fibers, subject to the condition that
$\sum_{i=1}^n m_i \langle c, F_i \rangle = \langle c_0, \varphi
\rangle$  for each connected component of each singular fiber.
Choosing these values to be positive, we apply Theorem~\ref{thm:5.2}(b)  to 
obtain symplectic forms  taming  $J$ (hence in $\cS$ as characterized 
in Theorem~\ref{contractible}(b)). 
Since $\langle f^*[\omega_\Sigma], F_i\rangle = 0$ for each
$F_i$, we obtain the required flexibility of $[\omega]$.
\end{proof} 

\begin{proof}[Proof of Theorem~\ref{contractible}] 
By Proposition~\ref{tame} and subsequent text, a structure $J$ satisfies the 
condition of Theorem~\ref{contractible}(b) if and only if it is 
$(\omega_\Sigma,f)$-compatible, such structures $J$ are easy to construct, 
and any $\omega$ as in (b) lies in $\cS$. 
Thus $\cS \ne \emptyset$ by Theorem~\ref{areas}. 
It remains to show that $\cS$ is contractible and any $\omega\in\cS$
tames some $J$ as in (b). 
As in the contractibility proof  of \cite[Theorem 2.11(b)]{G2}, 
a map of a sphere $S^m \to \cS$ can be
interpreted as a family $\omega$ of symplectic forms on the fibers of
the trivial bundle $S^m \times X$. 
Let $L \to S^m \times (X -K)$
denote the bundle of tangent spaces to fibers of $f$, and let
$L^\perp$ denote its orthogonal complement in $S^m \times T(X-K)$ with
respect to the given family $\omega$. 
(Then $L \oplus L^\perp = S^m
\times T(X-K)$ since $\omega$ is symplectic on each fiber of $f$, and $\omega$
is nondegenerate on each subbundle.)  
Let $J_0$ be a complex structure near $K$ induced by Definition~\ref{def:5.1}, 
pulled back to the $X$-fibers of $S^m\times X$.
Since $J_0 |S^m\times K = J^*$ is $\omega$-tame, $J_0$ is $\omega$-tame 
near some compact subset $\hat X \subset S^m\times X$ whose interior 
contains $S^m\times K$. 
Since the fibers of $f$ are $J_0$-holomorphic on $\hat X$,
we can obtain a new  complex structure $J_1$ on
$(S^m\times TX)|(\hat X-S^m\times K)$ by identifying $L^\perp$ there with 
$S^m\times TX/L$ and summing the
resulting complex structure on $L^\perp$ with the given one on $L$.
Then $J_1$ is $\omega$-tame by $\omega$-orthogonality of the line bundle 
$L^\perp$, and 
$f$ is both $J_0$- and $J_1$-holomorphic for the same complex 
structure on $S^m\times f^*T\Sigma |\hat X$. 
Thus, we can interpolate between $J_0$ and $J_1$  on $\hat X$ 
using Lemma~\ref{lem:3.2}(b):  
Set $T= \id_{S^m}\times df : S^m\times TX \to S^m
\times f^* T\Sigma$ restricted to $\hat X$.
Let each $W_x$ equal $\hat X$, with $J_0$ on $E|W_x$, and let 
$C= S^m\times K\cup \partial \hat X$, with $J_C$ given by $J_0$ near 
$S^m\times K$ and $J_1$ near $\partial\hat X$. 
We obtain an $\omega$-tame, $(\omega_\Sigma, f)$-compatible structure $J$ 
on $\hat X$ that agrees with $J_1$ on $\partial \hat X$. 
Since $L$ and $L^\perp$ are $J_1$-complex on $\partial\hat X$, 
we can extend $J$ over $S^m \times X$ by declaring $L$ and 
$L^\perp$ to be complex line bundles outside $\hat X$ as well. 
Replacing $S^m$ by a point in this argument gives the required $J$ on $X$ 
completing the proof of (b). 
(To smooth $J$ if desired, assume it agrees with the smooth $J_0$ near $K$, 
then perturb away from $K$ by first smoothing a $J$-complex line bundle 
complementary to $L$.) 
On the other hand, 
Lemma~\ref{lem:3.2}(a) (over $D^{m+1}\times X$ with $C = S^m \times X$)
extends $J$ to an $(\omega_\Sigma, f)$-compatible structure on
$D^{m+1} \times X$ that is $\omega$-tame on $\partial D^{m+1} \times X$. 
As in the uniqueness proof of Theorem~\ref{thm:5.2}(a) 
(which is the $m=0$ case), we can now construct a
$J$-taming symplectic form on each $\{p\} \times X$ 
(Theorem~\ref{areas}), then splice by
a partition of unity on $D^{m+1}$ to obtain a family of forms parametrized by
$D^{m+1}$. 
The resulting map $D^{m+1} \to \cS$ provides a nullhomotopy of the 
original map, showing $\pi_m(\cS)=0$ for all $m$. 
But $\cS$ is an open subset of a metrizable vector space of closed 
forms, so it is an ANR and hence contractible \cite{P} 
(cf.\ last paragraph of \cite[Section 3]{G2}).
\end{proof}

\begin{proof}[Proof of Lemma~\ref{lem:5.4}]
The lemma is obvious when $n=1$, so we prove the statement for fixed $n
> 1$, inductively assuming it for $n-1$.  
It is well-known that if
$G$ is a finite, connected graph for which each edge has two distinct
vertices, then there is a vertex $v$ with $G - \{v\}$ connected.  
(This is easy to prove by induction on the number of vertices: 
Fix a pair of vertices connected by an edge, delete all edges between them and
identify the pair, then apply the induction hypothesis to this smaller graph.)  
Given $A$ satisfying (a)--(c), reorder the coordinates so that
$G_A -\{v_n\}$ is connected. 
If (d) fails, then $\rank A < n$, but the
$(n-1) \times (n-1)$ matrix obtained from $A$ by deleting the last
row and column has $\rank n-1$ by induction.  
(Hypothesis (d) holds for
it since connectivity of $G_A$ implies $a_{nj} > 0$ for some $ j < n$.)  
Thus, we can assume $A$ satisfies (a)--(d).  
Note that $a_{nn} < 0$ by (a), (b) and (c). 
Eliminate the remaining entries in the last
column by row operations, adding $-a_{in}/a_{nn}$ times the $n^{th}$
row to the $i^{th}$ row for $i < n$, then delete the last row and
column and let $B = [b_{ij}]$ denote the resulting $(n-1) \times
(n-1)$ matrix.  
Since $b_{ij} = a_{ij} - \frac{a_{in}a_{nj}}{a_{nn}}$,
$B$ is symmetric, so it now suffices to show $B$ satisfies (a)--(d).
Condition (b) is obvious since $b_{ij} \geq a_{ij}$ for all $i,j <n$. 
Condition (a) follows since $G_B$ is obtained from $G_A$ by
deleting $v_n$ (and adjacent edges) and possibly adding edges.
Condition (c) follows from two applications of (c) for $A$.  
First we have (for $j = n$)
$$
\sum_{i=1}^{n-1} m_i a_{in} \leq -m_n a_{nn}.
$$
Since $a_{nn} < 0$, we obtain $\sum_{i=1}^{n-1} m_i b_{ij} =
\sum_{i=1}^{n-1} m_i a_{ij} - (\sum_{i=1}^{n-1} m_i a_{in})
\frac{a_{nj}}{a_{nn}} \leq \sum_{i=1}^{n-1} m_i a_{ij} + m_n a_{nj} =
\sum_{i=1}^n m_i a_{ij} \leq 0$ for all $j < n$, verifying (c).    
If the final inequality here is strict for some $j < n$, 
then (d) follows and we are finished.  
Otherwise (d) for $A$ implies that the displayed
inequality (for which $j = n$) is strict.  
By (a) there is some $j<n$ for which $a_{nj} \ne 0$, 
so the remaining inequality above is strict for this $j$, and (d) follows.  
(Note that condition (a) is crucial. 
Otherwise a matrix of diagonal blocks satisfying (b), (c) but
not (d) would be a counterexample.)
\end{proof}

\end{document}